\documentclass[12pt]{amsart}
\usepackage{amsmath,amsthm,amsfonts,amssymb,latexsym,enumerate} 
\usepackage{showlabels}

\newcommand{\CC}{{\mathbb{C}}}

\newcommand{\ZZ}{{\mathbb{Z}}}

\newcommand{\QQ}{{\mathbb{Q}}}

\newcommand{\SC} {\mathcal{S}}

\newcommand{\Tr}{{{\operatorname{Tr}}}}
\newcommand{\Gal}{{{\operatorname{Gal}}}}

\newcommand{\Irr}{{{\operatorname{Irr}}}}

\newcommand{\Ker}{\operatorname{Ker}}

\newtheorem*{thmA}{Theorem A}
\newtheorem*{conA'}{Conjecture A'}
\newtheorem*{thmB}{Theorem B}

\theoremstyle{definition}

\numberwithin{equation}{section}

\begin{document}

\title[A generalized character]{A generalized character\\ associated to element orders}

\author{Alexander Moret\'o}
\address{Departamento de Matem\'aticas, Universidad de Valencia, 46100
  Burjassot, Valencia, Spain}
\email{alexander.moreto@uv.es}

\thanks{I thank the anonymous referees for their careful reading of the paper. Research supported by  Ministerio de Ciencia e 
Innovaci\'on (Grants PID2019-103854GB-I00 and PID2022-137612NB-I00 funded by 
MCIN/AEI/10.13039/501100011033 and ``ERDF A way of making Europe")  and Generalitat Valenciana CIAICO/2021/163.}



\date{\today}

\begin{abstract}
Let $G$ be a finite group. We study the generalized character defined by  $\Xi(g)=|G|o(g)$, for $g\in G$, which is closely related to a function that has been very studied recently from a group theoretical point of view.
\end{abstract}

\maketitle


\section{Introduction}  

Let  $G$ be a finite group.  Given $g\in G$ we write $o(g)$ to denote the order of the element $g$.  The function 
$$
\psi(G)=\sum_{g\in G}o(g)
$$
was introduced in \cite{aji} and has been intensely studied since then. Its interest lies in the fact that it captures many structural properties of the group $G$. For instance, it was proved in \cite{aji} that if $G$ is a finite group of order $n$ then $\psi(G)\leq\psi(C_n)$. In other words, for any positive  integer $n$, the cyclic group of order $n$ is the  group of order $n$ whose  $\psi$ function is maximal among the groups of that order. See, for instance, \cite{hlm} and the references there for many other interesting properties of this function.

This has always been studied from a group theoretical point of view. Our goal in this note  is to introduce a character theoretical connection of the $\psi$ function. Our notation will follow \cite{isa}. We define the class function $$
\Xi(g)=|G|o(g)
$$ 
for $g\in G$. Sometimes, we will write $\Xi=\Xi_{(G)}$ if we want to make clear the group we are referring to. Note that $\psi(G)=[\Xi_{(G)},1_G]$.  We show the following.

\begin{thmA}
Let $G$ be a finite group. Then $|G|$ is the smallest positive integer $m=m(G)$ such that $mo(g)$ is a generalized character of $G$. Furthermore, $\Xi_{(G)}$ is a $\ZZ$-linear combination of characters $(1_C)^G$ for $C\leq G$ cyclic. 
\end{thmA}

We think that these properties of the order function are remarkable. 
 It is easy to see that if $G>1$ then $\Xi$ is never a character. Simply note that
$$
|G|=\Xi(1)=\sum_{\chi\in\Irr(G)}[\Xi,\chi]\chi(1)=\psi(G)+\sum_{\chi\in\Irr(G)-\{1_G\}}[\Xi,\chi]\chi(1)
$$
so $[\Xi,\chi]<0$ for some $\chi\in\Irr(G)$ (because $\psi(G)>|G|$). The generalized character $\Xi$ seems to have other interesting properties. Using GAP \cite{gap}, we have noticed that $\Xi$ tends to have many of the  irreducible characters of $G$ as irreducible constituents. Indeed, up to order 512, there is exactly one group $G$ with some $\chi\in\Irr(G)$ such that $[\Xi,\chi]=0$. It is ${\tt SmallGroup(504,8)}$. It may be worth remarking that the very much studied permutation character $\pi(g)=|C_G(g)|$ and its reciprocal class function $\tilde{\pi}(g)=|C_G(g)|^{-1}$ have a similar behavior that do not seem to be well understood yet. We refer the reader to \cite{tie, sam} and the references there.

We can compute in a rather explicit way the integer $[\Xi,\chi]$ for $\chi\in\Irr(G)$. If $\SC\subseteq G$, we write $\psi(\SC)=\sum_{x\in \SC}o(x)$. If $\chi\in\Irr(G)$ and $\alpha\in\CC$, we write $G_{\chi,\alpha}=\{g\in G\mid\chi(g)=\alpha\}$.

\begin{thmB}
Let $G$ be a finite group and let $\chi\in\Irr(G)$. Let $\{\alpha_1,\dots,\alpha_t\}$ be a complete set of representatives of the Galois-orbits of $\{\chi(x)\mid x\in G\}$. Then
$$
[\Xi,\chi]=\sum_{i=1}^t\psi(G_{\chi, \alpha_i})\Tr_{\QQ(\alpha_i)/\QQ}(\alpha_i).
$$
In particular, if $\lambda\in\Irr(G)$ is linear, $K=\Ker\lambda$, $m=|G:K|$, and $gK$  is a generator of $G/K$ (note that $G/K$ is cyclic), then $$[\Xi,\lambda]=\sum_{d\mid m}\psi(g^dK)\mu(m/d),$$
where $\mu$ is  the M\"obius function.
\end{thmB}

\section{Proofs}

Part of Theorem A is a consequence of a celebrated theorem of Artin.

\begin{proof}[Proof of Theorem A] 
The claim that $\Xi_{(G)}$ is a $\ZZ$-linear combination of characters $(1_C)^G$ for $C\leq G$ cyclic follows from Artin's theorem \cite{art}. 

Next, assume that $0<m<|G|$. We want to prove that $mo(g)$ is not a generalized character.  There exists a prime $p$ such that $m_p<|G|_p=p^a$. Write $m=p^bn$ with $(p,n)=1$, so that $b<a$. Consider the class function $$\mu(g)=mo(g)=\frac{p^bn}{|G|}\Xi(g)=\frac{n}{p^{a-b}|G|_{p'}}\Xi(g).$$
We want to see that $\mu$ is not a generalized character. Let $P$ be a Sylow $p$-subgroup of $G$. Then
$$
\mu_P(g)= \frac{n}{p^{a-b}|G|_{p'}}\Xi_P(g)=\frac{n}{p^{a-b}}\Xi_{(P)}(g).
$$
Now, note that $[\Xi_{(P)},1_P]=\psi(P)$ is a $p'$-number (congruent to $1$ modulo $p$). Hence
$$
[\mu_P,1_P]=\frac{n}{p^{a-b}}[\Xi_{(P)},1_P]
$$ 
is not an integer, so $\mu_P$ is not a generalized character. Hence $\mu$ is not a generalized character either, as desired.
\end{proof} 

Next, we prove  Theorem B.

\begin{proof}[Proof of Theorem B]
We have that
\begin{align*}
[\Xi,\chi]& =\sum_{x\in G}o(x)\chi(x)\\
& =\sum_{i=1}^t\sum_{\sigma\in\Gal(\QQ(\alpha_i)/\QQ)}\psi(G_{\chi, \alpha_i^{\sigma}})\alpha_i^{\sigma}\\
& =\sum_{i=1}^t\sum_{\sigma\in\Gal(\QQ(\alpha_i)/\QQ)}\psi(G_{\chi, \alpha_i})\alpha_i^{\sigma}\\
& =\sum_{i=1}^t\psi(G_{\chi, \alpha_i})\left(\sum_{\sigma\in\Gal(\QQ(\alpha_i)/\QQ)}\alpha_i^{\sigma}\right)\\
& =\sum_{i=1}^t\psi(G_{\chi, \alpha_i})\Tr_{\QQ(\alpha_i)/\QQ}(\alpha_i).
\end{align*}
where we have used in the third equality that the Galois group acts on the conjugacy classes of $G$ and Galois conjugation preserves the order of the elements.

Now, assume that $\lambda\in\Irr(G)$ is linear. Let $K=\Ker\lambda$, let $m=|G:K|$ and let $gK$ be a generator of $G/K$. 
Note that $\lambda(g)=\varepsilon$ for some  primitive $m$th root of unity $\varepsilon$. Furthermore,  the cosets $g^iK=\{g^iz\mid z\in K\}$ for $i=1,\dots,m$ form a partition of $G$ and $\lambda(g^iz)=\varepsilon^i$ for every $z\in K$.  We have that $\{\lambda(x)\mid x\in G\}=\{\varepsilon^i\mid i=1,\dots,m\}$ and a complete set of representatives of the Galois-orbits on this set is $\{\varepsilon^d\mid 1\leq d\leq m, d\mid m\}$. 
Hence
\begin{align*}
[\Xi,\lambda]& =\sum_{d\mid m}\psi(G_{\lambda, \varepsilon^d})\Tr_{\QQ(\varepsilon^d)/\QQ}(\varepsilon^d)\\
& =\sum_{d\mid m}\psi(g^dK)\mu(m/d),
\end{align*}
where in the last equality we have used that $G_{\lambda, \varepsilon^d}=g^dK$ and that if $\beta$ is a primitive $n$th root of unity then $\mu(n)=\Tr_{\QQ(\beta)/\QQ}(\beta)$.
\end{proof}







\begin{thebibliography}{99} 

\bibitem{aji}
H. Amiri, S. M. Jafarian Amiri, and I. M. Isaacs, Sums of element orders in finite groups, \textit{Comm. Algebra} \textbf{37} (2009), 2978--2980.

\bibitem{art}
E. Artin, Die gruppentheoretische Struktur der Diskriminanaten algebraischer Zahlk\"orper, \textit{J. reine angew. Math.} \textbf{164} (1931), 1--11.

\bibitem{hlm}
M. Herzog, P. Longobardi, M. Maj, Another criterion for solvability of finite groups, J. Algebra {\bf 597} (2022), 1--23.


\bibitem{gap}
The GAP group, \emph{{\sf GAP} - Groups, Algorithms, and Programming}.
  Version 4.4, 2004, {\sf http://www.gap-system.org}.

\bibitem{isa} I. M. Isaacs, \emph{Character Theory of Finite
  Groups}, Dover, New York, 1994. 
  
  
  \bibitem{sam}
  B. Sambale, The reciprocal character of the conjugation action, Publ. Math. Debrecen {\bf 99} (2021), 243--260.
  
  \bibitem{tie}
  P. H. Tiep, The conjugation representation of the binary modular congruence group,  (2023), arXiv:2303.02807.


 \end{thebibliography}
\end{document}